\documentclass{gtart}

\def\ifplaintex{\expandafter\ifx\csname documentclass\endcsname\relax}

\ifplaintex 
\else
\fi

\expandafter\ifx\csname beginpicture\endcsname\relax
\expandafter\ifx\csname documentclass\endcsname\relax
\input pictex \else
\input prepictex \input pictex \input postpictex \fi\fi

\def\gt{{\mathsurround=0pt\it $\cal G\mskip-2mu$eometry \&\ 
$\cal T\!\!$opology}}        

\def\gtp{{\mathsurround=0pt\it $\cal G\mskip-2mu$eometry \&\ 
$\cal T\!\!$opology $\cal P\!$ublications}}  

\def\lognumber#1{\def\thelognumber{#1}}
\def\volumenumber#1{\def\thevolumenumber{#1}}
\def\papernumber#1{\def\thepapernumber{#1}}
\def\volumeyear#1{\def\thevolumeyear{#1}}

\def\pagenumbers#1#2{\def\startpage{#1}\def\finishpage{#2}}
\def\published#1{\def\publishdate{#1}}
\def\proposed#1{\def\theproposer{#1}}
\def\seconded#1{\def\theseconders{#1}}
\def\received#1{\def\receiveddate{#1}}
\def\revised#1{\def\reviseddate{#1}}
\def\accepted#1{\def\accepteddate{#1}}
\def\asciititle#1{\def\theasciititle{#1}}

\def\coverauthors#1{\def\thecoverauthors{#1}}
\def\asciiauthors#1{\def\theasciiauthors{#1}}

\long\def\asciiabstract#1{\long\def\theasciiabstract{#1}}
\def\asciikeywords#1{\def\theasciikeywords{#1}}

\let\\\par\let\thelognumber\relax
\let\thevolumenumber\relax\let\thepapernumber\relax
\let\thevolumeyear\relax\let\thesamplenumber\relax\let\startpage\relax
\let\finishpage\relax\let\publishdate\relax\let\receiveddate\relax
\let\reviseddate\relax\let\accepteddate\relax\let\theasciititle\relax
\let\theasciiauthors\relax
\let\theasciiabstract\relax\let\theasciikeywords\relax
\let\theasciiemail\relax\let\theshortauthors\relax\let\theshorttitle\relax
\let\thecoverauthors\relax

\long\def\maketitlep{   

\count0=\startpage

\gt\hfill      
\beginpicture
\setcoordinatesystem units <0.33truein, 0.33truein> point at 2.2 0.9
\setplotsymbol ({$\cal G$})
\plotsymbolspacing=9truept
\circulararc 315 degrees from 0 1 center at 0 0
\setplotsymbol ({$\cal T$})
\circulararc 315 degrees from 1 -1 center at 1 0
\endpicture
\break
{\small\ifx\thesamplenumber\relax 
Volume \else Sample
\fi\thevolumenumber\ (\thevolumeyear)
\startpage--\finishpage\nl
Published: \publishdate}
\vglue 0.5truein plus 0.4fil minus 0.1truein

{\parskip=0pt\leftskip 0pt plus 1fil\def\\{\par\smallskip}{\ifplaintex\large
\else\Large\fi\bf\thetitle}\par\medskip}   

\vglue 0pt plus 0.1fil 

{\parskip=0pt\leftskip 0pt plus 1fil\def\\{\par}{\sc\theauthors}
\par\medskip}

\vglue 0pt plus 0.1fil 

{\small\parskip=0pt\let\newline\\
{\leftskip 0pt plus 1fil\def\\{\par}{\sl\theaddress}\par}
\expandafter\ifx\theemail\relax    
\relax\else\vglue 5pt plus 0.02fil minus 2pt\def\\{\stdspace{\rm 
and}\stdspace} 
\cl{Email:\stdspace\tt\theemail}\fi
\ifx\theurl\relax                  
\relax\else\vglue 5pt plus 0.02fil minus 2pt\def\\{\stdspace{\rm 
and}\stdspace}
\cl{URL:\stdspace\tt\theurl}\fi\par}

\vglue 7pt plus 0.3fil minus 3pt

{\bf Abstract}
\vglue 5pt plus 0.1fil minus 2pt

\theabstract

\vglue 7pt plus 0.3fil minus 3pt

{\bf AMS Classification numbers}\quad Primary:\quad \theprimaryclass

Secondary:\quad \thesecondaryclass

\vglue 5pt plus 0.3fil minus 2pt

{\bf Keywords}\quad \thekeywords

\vglue 10pt plus 0.5fil minus 5pt

{\small  Proposed: \theproposer\hfill Received: \receiveddate\nl
Seconded: \theseconders\hfill 
\ifx\reviseddate\relax                         
Accepted: \accepteddate                        
\else
Revised: \reviseddate                          
\fi}
\eject
}       

\let\maketitlepage\maketitlep
\let\maketitle\maketitlepage

\font\phead=cmsl9 scaled 950
\font=cmsl9 scaled 1050
\font\pnum=cmbx10 scaled 913
\font\lnum=cmbx10 
\font\pfoot=cmsl9 scaled 950
\font=cmsl9 scaled 1050
\ifplaintex
\headline{\vbox to 0pt{\vskip -4.5mm\line{\small\phead\ifnum
\count0=\startpage ISSN 1364-0380 (on line)
1465-3060 (printed) \hfill {\pnum\folio}\else\ifodd\count0\def\\{ }%
\ifx\theshorttitle\relax\thetitle\else\theshorttitle\fi\hfill{\pnum\folio}
\else\def\\{ and }{\pnum\folio}\hfill\ifx\theshortauthors\relax\theauthors
\else\theshortauthors\fi\fi\fi}\vss}}
\footline{\vbox to 0pt{\vglue 0mm\line{\small\pfoot\ifnum\count0=\startpage
\copyright\ \gtp\hfill\else
\gt, Volume \thevolumenumber\ (\thevolumeyear)\hfill\fi}\vss
}}
\else
\makeatletter
\def\@oddhead{{\small\lhead\ifnum\count0=\startpage ISSN 1364-0380 (on line)
1465-3060 (printed) \hfill {\lnum\number\count0}\else\ifodd\count0
\def\\{ }\ifx\theshorttitle\relax \thetitle \else\theshorttitle\fi\hfill
{\lnum\number\count0}\else\def\\{ and }{\lnum\number\count0}
\hfill\ifx\theshortauthors\relax 
\theauthors\else\theshortauthors\fi\fi\fi}}\def\@evenhead{\@oddhead}
\def\@oddfoot{\small\lfoot\ifnum\count0=\startpage\copyright\ \gtp\hfill\else
\gt, Volume \thevolumenumber\ (\thevolumeyear)\hfill\fi}
\def\@evenfoot{\@oddfoot}
\makeatother
\fi

\newwrite\gtoutfile
\long\gdef\makeheadfile{  
{\def\\{, }\def\s{ }
\immediate\openout\gtoutfile head.xxx
\immediate\write\gtoutfile{To: math@arxiv.org}
\immediate\write\gtoutfile{Subject: put or rep NNNNN:pppp}
\immediate\write\gtoutfile{--text follows this line--}
\immediate\write\gtoutfile{Proxy-for: \ifx\theasciiauthors\relax
\theauthors\else\theasciiauthors\fi\s<\ifx\theasciiemail\relax\theemail\else\theasciiemail\fi>}
\immediate\write\gtoutfile{\noexpand\\}
\immediate\write\gtoutfile{Authors: \ifx\theasciiauthors\relax
\theauthors\else\theasciiauthors\fi}
{\def\\{ }\immediate\write\gtoutfile{Title: \ifx\theasciititle\relax
\thetitle\else\theasciititle\fi}}
\immediate\write\gtoutfile{Subj-class: GT or SG or MG etc}
\immediate\write\gtoutfile{MSC-class: \theprimaryclass\ifx\thesecondaryclass\relax\else, \thesecondaryclass\fi}
\immediate\write\gtoutfile{Journal-ref: Geom. Topol. \thevolumenumber
(\thevolumeyear) \startpage-\finishpage}
\immediate\write\gtoutfile{Comments: Published by Geometry and Topology at}
\immediate\write\gtoutfile{\s\s http://www.maths.warwick.ac.uk/gt/GTVol\thevolumenumber/paper\thepapernumber.abs.html}
\immediate\write\gtoutfile{\noexpand\\}
\immediate\write\gtoutfile{}
\ifx\theasciiabstract\relax
\immediate\write\gtoutfile{\theabstract}\else
\immediate\write\gtoutfile{\theasciiabstract}\fi
\immediate\write\gtoutfile{}
\immediate\write\gtoutfile{\noexpand\\}
\immediate\write\gtoutfile{}
\immediate\closeout\gtoutfile}}  

\def\maketitlepage{\maketitlep\makeheadfile}
\let\maketitle\maketitlepage

\def\ifplaintex{\expandafter\ifx\csname documentclass\endcsname\relax}

\ifplaintex 
\else
\fi

\expandafter\ifx\csname beginpicture\endcsname\relax
\expandafter\ifx\csname documentclass\endcsname\relax
\input pictex \else
\input prepictex \input pictex \input postpictex \fi\fi

\def\gt{{\mathsurround=0pt\it $\cal G\mskip-2mu$eometry \&\ 
$\cal T\!\!$opology}}        

\def\gtp{{\mathsurround=0pt\it $\cal G\mskip-2mu$eometry \&\ 
$\cal T\!\!$opology $\cal P\!$ublications}}  

\def\lognumber#1{\def\thelognumber{#1}}
\def\volumenumber#1{\def\thevolumenumber{#1}}
\def\papernumber#1{\def\thepapernumber{#1}}
\def\volumeyear#1{\def\thevolumeyear{#1}}

\def\pagenumbers#1#2{\def\startpage{#1}\def\finishpage{#2}}
\def\published#1{\def\publishdate{#1}}
\def\proposed#1{\def\theproposer{#1}}
\def\seconded#1{\def\theseconders{#1}}
\def\received#1{\def\receiveddate{#1}}
\def\revised#1{\def\reviseddate{#1}}
\def\accepted#1{\def\accepteddate{#1}}
\def\asciititle#1{\def\theasciititle{#1}}

\def\coverauthors#1{\def\thecoverauthors{#1}}
\def\asciiauthors#1{\def\theasciiauthors{#1}}

\long\def\asciiabstract#1{\long\def\theasciiabstract{#1}}
\def\asciikeywords#1{\def\theasciikeywords{#1}}

\let\\\par\let\thelognumber\relax
\let\thevolumenumber\relax\let\thepapernumber\relax
\let\thevolumeyear\relax\let\thesamplenumber\relax\let\startpage\relax
\let\finishpage\relax\let\publishdate\relax\let\receiveddate\relax
\let\reviseddate\relax\let\accepteddate\relax\let\theasciititle\relax
\let\theasciiauthors\relax
\let\theasciiabstract\relax\let\theasciikeywords\relax
\let\theasciiemail\relax\let\theshortauthors\relax\let\theshorttitle\relax
\let\thecoverauthors\relax

\long\def\maketitlep{   

\count0=\startpage

\gt\hfill      
\beginpicture
\setcoordinatesystem units <0.33truein, 0.33truein> point at 2.2 0.9
\setplotsymbol ({$\cal G$})
\plotsymbolspacing=9truept
\circulararc 315 degrees from 0 1 center at 0 0
\setplotsymbol ({$\cal T$})
\circulararc 315 degrees from 1 -1 center at 1 0
\endpicture
\break
{\small\ifx\thesamplenumber\relax 
Volume \else Sample
\fi\thevolumenumber\ (\thevolumeyear)
\startpage--\finishpage\nl
Published: \publishdate}
\vglue 0.5truein plus 0.4fil minus 0.1truein

{\parskip=0pt\leftskip 0pt plus 1fil\def\\{\par\smallskip}{\ifplaintex\large
\else\Large\fi\bf\thetitle}\par\medskip}   

\vglue 0pt plus 0.1fil 

{\parskip=0pt\leftskip 0pt plus 1fil\def\\{\par}{\sc\theauthors}
\par\medskip}

\vglue 0pt plus 0.1fil 

{\small\parskip=0pt\let\newline\\
{\leftskip 0pt plus 1fil\def\\{\par}{\sl\theaddress}\par}
\expandafter\ifx\theemail\relax    
\relax\else\vglue 5pt plus 0.02fil minus 2pt\def\\{\stdspace{\rm 
and}\stdspace} 
\cl{Email:\stdspace\tt\theemail}\fi
\ifx\theurl\relax                  
\relax\else\vglue 5pt plus 0.02fil minus 2pt\def\\{\stdspace{\rm 
and}\stdspace}
\cl{URL:\stdspace\tt\theurl}\fi\par}

\vglue 7pt plus 0.3fil minus 3pt

{\bf Abstract}
\vglue 5pt plus 0.1fil minus 2pt

\theabstract

\vglue 7pt plus 0.3fil minus 3pt

{\bf AMS Classification numbers}\quad Primary:\quad \theprimaryclass

Secondary:\quad \thesecondaryclass

\vglue 5pt plus 0.3fil minus 2pt

{\bf Keywords}\quad \thekeywords

\vglue 10pt plus 0.5fil minus 5pt

{\small  Proposed: \theproposer\hfill Received: \receiveddate\nl
Seconded: \theseconders\hfill 
\ifx\reviseddate\relax                         
Accepted: \accepteddate                        
\else
Revised: \reviseddate                          
\fi}
\eject
}       

\let\maketitlepage\maketitlep
\let\maketitle\maketitlepage

\font\phead=cmsl9 scaled 950
\font=cmsl9 scaled 1050
\font\pnum=cmbx10 scaled 913
\font\lnum=cmbx10 
\font\pfoot=cmsl9 scaled 950
\font=cmsl9 scaled 1050
\ifplaintex
\headline{\vbox to 0pt{\vskip -4.5mm\line{\small\phead\ifnum
\count0=\startpage ISSN 1364-0380 (on line)
1465-3060 (printed) \hfill {\pnum\folio}\else\ifodd\count0\def\\{ }%
\ifx\theshorttitle\relax\thetitle\else\theshorttitle\fi\hfill{\pnum\folio}
\else\def\\{ and }{\pnum\folio}\hfill\ifx\theshortauthors\relax\theauthors
\else\theshortauthors\fi\fi\fi}\vss}}
\footline{\vbox to 0pt{\vglue 0mm\line{\small\pfoot\ifnum\count0=\startpage
\copyright\ \gtp\hfill\else
\gt, Volume \thevolumenumber\ (\thevolumeyear)\hfill\fi}\vss
}}
\else
\makeatletter
\def\@oddhead{{\small\lhead\ifnum\count0=\startpage ISSN 1364-0380 (on line)
1465-3060 (printed) \hfill {\lnum\number\count0}\else\ifodd\count0
\def\\{ }\ifx\theshorttitle\relax \thetitle \else\theshorttitle\fi\hfill
{\lnum\number\count0}\else\def\\{ and }{\lnum\number\count0}
\hfill\ifx\theshortauthors\relax 
\theauthors\else\theshortauthors\fi\fi\fi}}\def\@evenhead{\@oddhead}
\def\@oddfoot{\small\lfoot\ifnum\count0=\startpage\copyright\ \gtp\hfill\else
\gt, Volume \thevolumenumber\ (\thevolumeyear)\hfill\fi}
\def\@evenfoot{\@oddfoot}
\makeatother
\fi

\newwrite\gtoutfile
\long\gdef\makeheadfile{  
{\def\\{, }\def\s{ }
\immediate\openout\gtoutfile head.xxx
\immediate\write\gtoutfile{To: math@arxiv.org}
\immediate\write\gtoutfile{Subject: put or rep NNNNN:pppp}
\immediate\write\gtoutfile{--text follows this line--}
\immediate\write\gtoutfile{Proxy-for: \ifx\theasciiauthors\relax
\theauthors\else\theasciiauthors\fi\s<\ifx\theasciiemail\relax\theemail\else\theasciiemail\fi>}
\immediate\write\gtoutfile{\noexpand\\}
\immediate\write\gtoutfile{Authors: \ifx\theasciiauthors\relax
\theauthors\else\theasciiauthors\fi}
{\def\\{ }\immediate\write\gtoutfile{Title: \ifx\theasciititle\relax
\thetitle\else\theasciititle\fi}}
\immediate\write\gtoutfile{Subj-class: GT or SG or MG etc}
\immediate\write\gtoutfile{MSC-class: \theprimaryclass\ifx\thesecondaryclass\relax\else, \thesecondaryclass\fi}
\immediate\write\gtoutfile{Journal-ref: Geom. Topol. \thevolumenumber
(\thevolumeyear) \startpage-\finishpage}
\immediate\write\gtoutfile{Comments: Published by Geometry and Topology at}
\immediate\write\gtoutfile{\s\s http://www.maths.warwick.ac.uk/gt/GTVol\thevolumenumber/paper\thepapernumber.abs.html}
\immediate\write\gtoutfile{\noexpand\\}
\immediate\write\gtoutfile{}
\ifx\theasciiabstract\relax
\immediate\write\gtoutfile{\theabstract}\else
\immediate\write\gtoutfile{\theasciiabstract}\fi
\immediate\write\gtoutfile{}
\immediate\write\gtoutfile{\noexpand\\}
\immediate\write\gtoutfile{}
\immediate\closeout\gtoutfile}}  

\def\maketitlepage{\maketitlep\makeheadfile}
\let\maketitle\maketitlepage

\lognumber{252}
\volumenumber{6}
\papernumber{11}
\volumeyear{2002}
\pagenumbers{355}{360}
\received{5 May 2002}
\revised{19 June 2002}
\accepted{26 June 2002}
\published{27 June 2002}

\proposed{Robion Kirby}
\seconded{Walter Neumann, Vaughan Jones}

\input{rlepsf}
\let\relabela\adjustrelabel

\newtheorem{theorem}{Theorem}
\newtheorem{lemma}[theorem]{Lemma}
\newtheorem{proposition}[theorem]{Proposition}
\newtheorem{conjecture}[theorem]{Conjecture}

\theoremstyle{definition}
\newtheorem{remark}[theorem]{Remark}
\newtheorem{definition}[theorem]{Definition}

\nocolon

\begin{document}

\title{Burnside obstructions to the Montesinos--Nakanishi\\3--move 
conjecture}
\asciititle{Burnside obstructions to the Montesinos-Nakanishi\\3-move 
conjecture}
\authors{Mieczys{\l}aw K~D{\c a}bkowski\\J\'ozef H~Przytycki}
\coverauthors{Mieczys{\noexpand\l}aw K D{\noexpand\c a}bkowski\\J\noexpand\'ozef H Przytycki}
\asciiauthors{Mieczyslaw K Dabkowski\\Jozef H Przytycki}

\address{Department of Mathematics, The George Washington University\\
Washington, DC 20052, USA}
\email{mdab@gwu.edu, przytyck@gwu.edu}

\begin{abstract}
Yasutaka Nakanishi asked in 1981 whether a 3--move is an unknotting
operation. In Kirby's problem list, this question is called
{\it The Montesinos--Nakanishi 3--move conjecture}. We define
the $n$th Burnside group of a link and use the 3rd
Burnside group to answer Nakanishi's question; ie, we show that some
links cannot be reduced to trivial links by 3--moves.
\end{abstract}

\asciiabstract{%
Yasutaka Nakanishi asked in 1981 whether a 3-move is an unknotting
operation. In Kirby's problem list, this question is called
`The Montesinos-Nakanishi 3-move conjecture'. We define
the n-th Burnside group of a link and use the 3rd
Burnside group to answer Nakanishi's question; ie, we show that some
links cannot be reduced to trivial links by 3-moves.}

\keywords{Knot, link, tangle, 3--move, rational move, braid, Fox
coloring, Burnside group, Borromean rings, Montesinos--Nakanishi
conjecture, branched cover, core group, lower central series,
associated graded Lie ring, skein module}

\asciikeywords{Knot, link, tangle, 3-move, rational move, braid, Fox
coloring, Burnside group, Borromean rings, Montesinos-Nakanishi
conjecture, branched cover, core group, lower central series,
associated graded Lie ring, skein module}

\primaryclass{57M27}\secondaryclass{20D99}

\maketitlepage

One of the oldest elementary formulated problems in classical Knot Theory is 
the 3--move conjecture of Nakanishi. A 3--move on a link is a local change that 
involves replacing parallel lines by 3 half-twists (Figure 1).

\begin{figure}[ht!]
\epsfysize 1.4cm
\centerline{\epsfbox{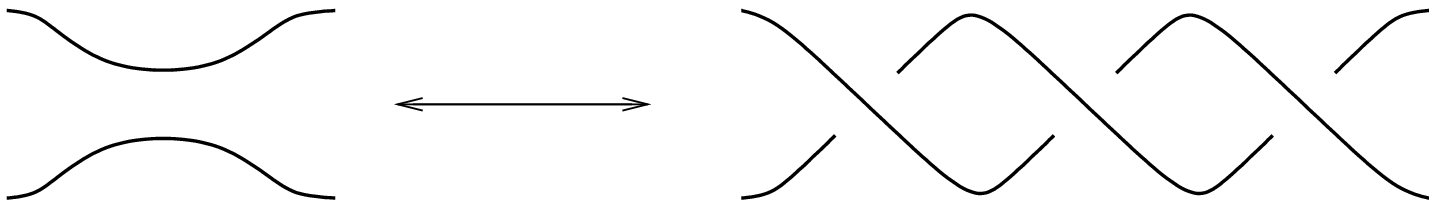}}
\caption{}
\end{figure}

\begin{conjecture}
{\rm(Montesinos--Nakanishi, Kirby's problem
list; Problem 1.59(1), \cite{Kir})}\label{1}\qua
Any link can be reduced to a trivial link by a sequence of 3--moves.
\end{conjecture}

The conjecture has been proved to be valid for several classes of links
by Chen, 
Nakanishi, Przytycki and Tsukamoto (eg, closed 4--braids and 4--bridge 
links).

Nakanishi, in 1994, and Chen, in 1999, have presented examples of links which 
they were not able to reduce: $L_{2BR}$, the 2--parallel of the Borromean 
rings, and $\hat\gamma$,
the closure of the square of the center of the fifth braid group,
ie, $\gamma = (\sigma_1 \sigma_2 \sigma_3 \sigma_4)^{10}$.

\begin{remark}\label{Remark 2}
In \textup{\cite{Pr-1}} it was noted that 3--moves preserve the first
homology of the double branched cover of a link $L$ with $Z_3$
coefficients ($H_1(M_L^{(2)};Z_3)$). Suppose that $\hat\gamma$
(respectively $L_{2BR}$) can be reduced by 3--moves to the trivial link
$T_n$. Since $H_1(M_{\hat\gamma}^{(2)};Z_3) = Z_3^4$,
$H_1(M_{L_{2BR}}^{(2)};Z_3) = Z_3^5$ and $H_1(M_{T_n}^{(2)};Z_3) =
Z_3^{n-1}$ where $T_n$ is a trivial link of $n$ components, it follows
that $n=5$ (respectively n=6).
\end{remark}

We show below that neither $\hat\gamma$ nor $L_{2BR}$ can be reduced
by 3--moves to trivial links.

The tool we use is a non-abelian version of Fox $n$--colorings, which we
shall call the $n$th Burnside group of a link, $B_L(n)$.

\begin{definition}\label{Def. 3} The $n$th Burnside group of a link is the 
quotient of the
fundamental group of the double branched cover of $S^3$ with the link as
the branch set divided
by all relations of the form $a^n=1$. Succinctly:\ 
$B_L(n)=\pi_1(M_L^{(2)})/(a^n)$.
\end{definition}

\begin{proposition}\label{Prop.4}
$B_L(3)$ is preserved by 3--moves.
\end{proposition}

\begin{proof} In the proof we use the core group interpretation of 
$\pi_1(M_L^{(2)})$. Let $D$ be a diagram of a link $L$. We define (after 
\cite{Joy,F-R}) the associated core group $\Pi^{(2)}_D$ of $D$ as follows: 
generators of $\Pi^{(2)}_D$
correspond to arcs of the diagram.
Any crossing $v_s$ yields
the relation $r_s=y_iy_j^{-1}y_iy_k^{-1}$ where $y_i$ corresponds
to the overcrossing
and $y_j,y_k$ correspond to the undercrossings at $v_s$ (see Figure 2).
In this presentation of $\Pi^{(2)}_L$ one relation can be dropped since it is 
a consequence of others.
Wada proved that $\Pi^{(2)}_D = \pi_1(M_L^{(2)})\ast Z$, \cite{Wa} 
(see \cite{Pr-2} for an elementary proof
using only Wirtinger presentation).
Furthermore, if we put $y_i=1$ for any fixed generator, then
$\Pi^{(2)}_D$ reduces to $\pi_1(M_L^{(2)})$. The last part of our proof is 
illustrated in Figure 2. \end{proof}

\begin{figure}[ht!]
\relabelbox\small
\epsfysize 3cm
\centerline{\epsfbox{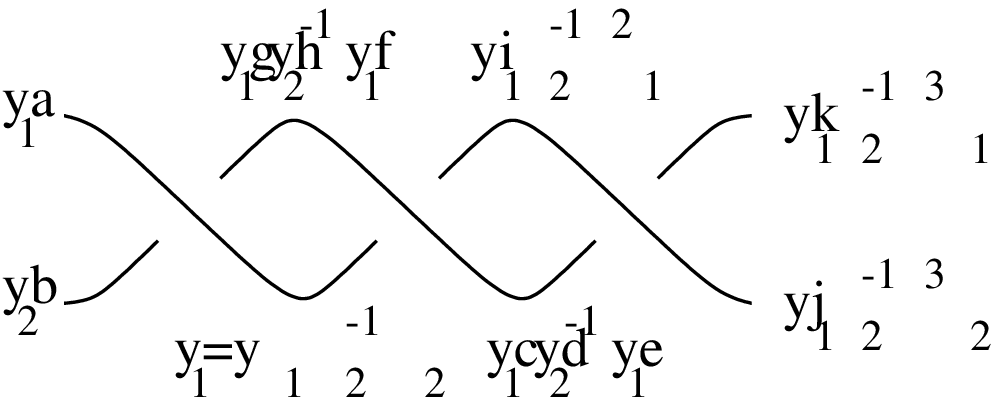}}
\relabela <-3pt,0pt> {y=y}{$y_1{=}y_1y_2^{-1}y_2$}
\relabela <-3pt,0pt> {ya}{$y_1$}
\relabela <-3pt,0pt> {yb}{$y_2$}
\relabela <-3pt,0pt> {yk}{$(y_1y_2^{-1})^3y_1$}
\relabela <-12pt,0pt> {yd}{$y_1y_2^{-1}y_1$}
\relabela <-3pt,0pt> {yg}{$y_1y_2^{-1}y_1$}
\relabela <-3pt,0pt> {yi}{$(y_1y_2^{-1})^2y_1$}
\relabela <-3pt,0pt> {yj}{$(y_1y_2^{-1})^3y_2$}
\endrelabelbox
\caption{}
\end{figure}

\begin{lemma}\label{Lemma 5}
$B_{\hat\gamma}(3) = \{x_1,x_2,x_3,x_4 \mid a^3 \mbox{ for any word }a, \
P_1, P_2, P_3, P_4 \}$, where
$$P_i=
x_1x_2^{-1}x_3x_4^{-1}x_1^{-1}x_2x_3^{-1}x_4x_i
x_4x_3^{-1}x_2x_1^{-1}x_4^{-1}x_3x_2^{-1}x_1x_i^{-1} .$$
\end{lemma}
\begin{proof}
Consider the 5--braid $\gamma = (\sigma_1 \sigma_2 \sigma_3 \sigma_4)^{10}$
(Figure 3). If we label initial arcs of the braid by $x_1,x_2,x_3,x_4$ and
$x_5$, and use core relations (progressing from left to right) we
obtain labels $Q_1,Q_2,Q_3,Q_4$ and $Q_5$ on the final arcs of the
braid where
$$Q_i = x_1x_2^{-1}x_3x_4^{-1}x_5x_1^{-1}x_2x_3^{-1}x_4x_5^{-1}
x_i
x_5^{-1}x_4x_3^{-1}x_2x_1^{-1}x_5x_4^{-1}x_3x_2^{-1}x_1 .$$
For a group $\Pi^{(2)}_{\hat\gamma}$, of the closed braid $\hat\gamma$,
we have relations $Q_i=x_i$. To obtain $\pi_1(M^{(2)}_{\hat\gamma})$ we
can put $x_5=1$, and delete one relation, say $Q_5x_5^{-1}$.
These lead to the presentation of $B_{\hat\gamma}(3)$ described in the lemma.
\end{proof}

\begin{figure}[ht!]
\relabelbox\small
\epsfysize 1.7cm
\centerline{\epsfbox{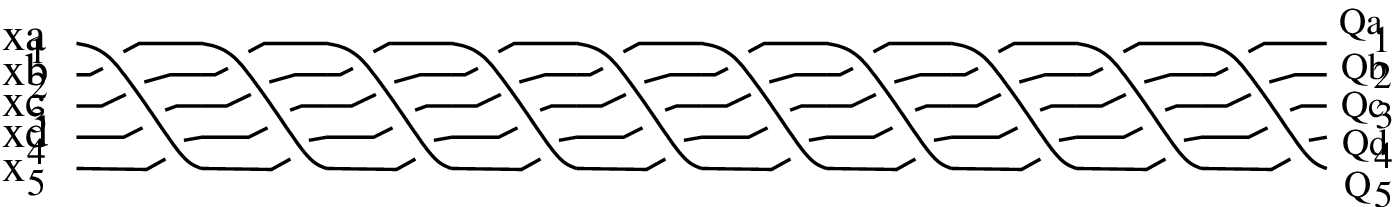}}
\relabel{x}{\scriptsize$x_5$}
\relabel{xa}{\scriptsize$x_1$}
\relabel{xb}{\scriptsize$x_2$}
\relabel{xc}{\scriptsize$x_3$}
\relabel{xd}{\scriptsize$x_4$}
\relabel{Q}{\scriptsize$Q_5$}
\relabel{Qa}{\scriptsize$Q_1$}
\relabel{Qb}{\scriptsize$Q_2$}
\relabel{Qc}{\scriptsize$Q_3$}
\relabel{Qd}{\scriptsize$Q_4$}
\endrelabelbox
\caption{}
\end{figure}

\begin{theorem}\label{Theorem 6}
The links $\hat\gamma$ and $L_{2BR}$ are not 3--move reducible to trivial 
links.
\end{theorem}
\begin{proof}
Let $B(n,3)$ denote the classical free $n$ generator Burnside group 
of exponent $3$. 
As shown by Burnside \cite{Bu},
$B(n,3)$ is a finite group. Its order, $|B(n,3)|$, is equal to
$3^{n + {n\choose 2} + {n\choose 3}}$.
For a trivial link: $B_{T_k}(3)= B(k-1,3)$. In order to prove that
$\hat\gamma$ and $L_{2BR}$ are not 3--move reducible to trivial links,
it suffices to show that $B_{\hat\gamma}(3) \neq B(4,3)$ and
$B_{L_{2BR}}(3) \neq B(5,3)$ (see Remark 2). We have demonstrated these 
to be true both by manual computation, and by 
using the programs GAP, Magnus and Magma. 
More details in the case of $\hat\gamma$
are provided below.

For the manual calculations, one first observes that
for any $i$, $P_i$ is in the third term of the lower central series
of $B(4,3)$. In particular, for $u=x_1x_2^{-1}x_3x_4^{-1}$ and $\bar u 
=x_1^{-1}x_2x_3^{-1}x_4$, one has $u\bar{u} \in [B(4,3),B(4,3)]$
and $P_i = [u\bar u,x_i\bar u]$. It is known (\cite{V-Lee}),
that $B(4,3)$ is of class
3 (the lower central series has 3 terms), and that the third term is
isomorphic to $Z_3^4$ with basis: $e_1=[[x_2,x_3],x_4]$,
$e_2=[[x_1,x_3],x_4]$, $e_3=[[x_1,x_2],x_4]$ and $e_4=[[x_1,x_2],x_3]$.
It now takes an elementary linear algebra calculation (see Lemma 7 below)
 to show that $P_1, P_2, P_3, P_4$
form another basis of the third term of the lower central series of $B(4,3)$.
Thus $|B_{\hat\gamma}(3)| = 3^{10}$.
\end{proof}
\begin{lemma}\label{Lemma 7}
$P_1,P_2,P_3$, and $P_4$ form a basis of the third term of the lower central
series of $B(4,3)$.
\end{lemma}

\begin{proof}
In the associated graded Lie ring $L(4,3)$ of $B(4,3)$ (\cite{V-Lee}),
the third term (denoted $L_3$) is isomorphic to $Z_3^4$ with basis
$e_1,e_2,e_3,e_4$. In $L(4,3)$, which is a linear space over $Z_3$,
one uses an additive notation and the bracket in the group becomes
a (non-associative) product (\cite{V-Lee}). In this notation
$e_1=x_2x_3x_4$, $e_2=x_1x_3x_4$, $e_3=x_1x_2x_4$ and $e_4=x_1x_2x_3$.
In the calculation expressing $P_i$ in the basis we use the
following identities in $L_3$ (\cite{V-Lee}; page 89).
$$xyzt=0, xyz=yzx=zxy=-xzy=-zyx=-yxz, xyy=0.$$
Now we have: $P_i= (u{\bar u})(x_i{\bar u})(u{\bar u})^{-1}(x_i{\bar u})^{-1}=
[(u{\bar u})^{-1},(x_i{\bar u})^{-1}] = [u{\bar u},x_i{\bar u}]$ as
the last term of the lower central series is in the center of $B(4,3)$.
Furthermore, we have $u{\bar u}= x_1x_2^{-1}x_3x_4^{-1}x_1^{-1}x_2x_3^{-1}x_4 =
[x_2^{-1}x_3x_4^{-1},x_1^{-1}][x_3x_4^{-1},x_2][x_4^{-1},x_3^{-1}]$.
Writing $P_i$ additively in $L_3$ one obtains:
$$P_i=((-x_2+x_3-x_4)(-x_1) + (x_3-x_4)x_2 + x_4x_3)(x_i-x_1+x_2-x_3+x_4).$$
After simplifications one gets:
$$P_1= -e_1, P_2= e_1 +e_2,  P_3= e_1 -e_2 -e_3,\ and \ P_4=e_1 -e_2 +
e_3+e_4.$$  
The matrix expressing $P_i$'s in terms of $e_i$'s is the upper triangular
matrix with the determinant equal to $1$.  Therefore the lemma follows.
\end{proof}

A similar calculation establishes that $|B_{L_{2BR}}(3)| < |B(5,3)|$.
$B(5,3)$ is of class 3 and has $3^{25}$ elements.
Considering $L_{2BR}$ as a closed 6--braid 
we note that $B_{L_{2BR}}(3)$ is obtained from $B(5,3)$ by adding 
5 relations $R_1,...,R_5$. 
Relations $\{R_i\}$ are in the last term of the lower central series of $B(5,3)$
(and of the associated graded algebra $L(5,3))$. Relations form a 
4--dimensional 
subspace in $L_3=Z_3^{10}$.  Thus $|B_{L_{2BR}}(3)|= 3^{21}$.

For a computer verification showing that $B_{\hat\gamma}(3) \neq B(4,3)$ consider
any presentation of $B(4,3)$ (eg, Magma solution by Mike Newman \cite{New})
and add the relations $P_i$ to obtain a presentation
of $B_{\hat\gamma}(3)$. Using any of the algebra programs mentioned above,
one verifies that $|B_{\hat\gamma}(3)| = 3^{10}$ while
$|B(4,3)|=3^{14}$.

The solution of the Nakanishi--Montesinos 3--move conjecture, presented
above, is the first instance of application of Burnside groups of links. 
It was motivated by the analysis of cubic skein modules of 3--manifolds.
The next step is the application of Burnside groups to rational 
moves on links.  
This, in turn, should have deep implications to the 
theory of skein modules \cite{Pr-2}.

\end{document}